# Helices associated to helical curves, relatively normal-slant helices and isophote curves

**Mehmet Önder**
*Delibekirli Village, Tepe Street, No: 63, 31440, Kırıkhan, Hatay, Turkey.*
*E-mail: mehmetonder197999@gmail.com*
*Orcid Id: https://orcid.org/0000-0002-9354-5530*

**Abstract**
This study introduces a new type of general helix called associated helix which is associated to a special surface curve. The basic idea is to determinate the parametric form of an associated helix by means of Darboux frame and surface curvatures of a special surface curve such as helical curve, relatively normal-slant helix or isophote curve. For each surface curve, a differential equation system is obtained and by solving this system, parametric form of an associated helix is introduced.

## 1. Introduction

Helix curve(or general helix) is the most common curve type in science and nature. In nature, these curves are frequently encountered in subjects such as DNA helixes, stem structures of plants like beans and ivy. Moreover, In many applications, helices are used [3,5,7,21]. For instance, these curves can be used to simulate the kinematic motion, design the highways and describe tool paths in CAD and computer graphics [24]. Then, helices and their applications are studied widely [1,3,6,11,18].

A helix is defined by property that its tangent always makes a constant angle with a constant direction called axis of the helix and characterized by Lancret's theorem saying that the function $\tau/\kappa$ is constant, where $\kappa, \tau$ are first and second curvatures of curve, respectively[4,23].

In the definition of general helix, considering unit principal normal vector instead of tangent vector, a new special curve is obtained which is called slant helix. Slant helix was defined by Izumiya and Takeuchi in 2004 [13]. Later, Kula and Yaylı has introduced spherical indicatrix of a slant helix [15]. Ali has given position vector of a slant helix by means of its first curvature function $\kappa(s)$ [2]. Furthermore, by taking the Darboux vector instead of tangent vector, Zıplar, Şenol and Yaylı have defined Darboux helices and showed that a curve is a Darboux helix if and only if it is a slant helix[25].

Moreover, a curve $\alpha$ lying fully on an oriented surface is called a surface curve and has an orthonormal frame $\{T, V, U\}$ along the curve. This new frame is called Darboux frame, where $T$ and $U$ are unit tangent vector and unit surface normal along $\alpha$, respectively, and $V = T \times U$ is a unit vector field along the curve $\alpha$. Such a surface curve is called helical curve(relatively normal-slant helix or isophote curve, respectively) if the vector field $T$ (vector filed $V$ or vector field $U$, respectively) makes a constant angle with a constant(or fixed) direction. Öztürk and Hacısalihoğlu defined such curves as strip slant helices[19]. Helical surface curve has been studied by Puig Pey, Galvez and Iglesias and they have introduced a new method for generating the helical tool-paths for both parametric form and implicit form of a surface[21]. The definition and characterizations of relatively normal-slant helix have been introduced by Macit and Düldül and they have also given the axis of this curve[16].

Furthermore, the geometric description of an isophote curve has been introduced by Doğan and Yaylı[10]. They have also studied isophote curves in Minkowski space[8,9]. Just like helices, these curves have a wide range of applications. In the optics branch of physics, an isophote curve on a surface is a nice consequence of Lambert's cosine law stating that the intensity of illumination on a diffuse surface is proportional to the cosine of the angle generated between the surface normal vector $U$ and the light vector $d$. From this law, it is clear that the intensity is irrespective of the actual viewpoint and it means the illumination is the same when viewed from any direction [10]. Moreover, the points of isophote curves on a surface have the same light intensity from a given source. By detecting irregularities along these curves on a free form surface, isophote curves are used in car body construction[20]. Sara researched local shading of a surface through isophote curve properties and studied on accurate estimation of surface normal tilt and on qualitatively correct Gaussian curvature recovery and for this purpose, he used the fundamental theory of surfaces[22]. Furthermore, the parametrizations of isophote curves on a canal surface or a rotation surface have been given by Kim and Lee[14].

In this study, we define new types of associated helices. We name these new helices as surface curve-connected(SCC) associated helices and we determinate parametrizations of SCC-associated helices by considering helical curves, relatively normal-slant helices and isophote curves on an oriented surface.

## 2. Preliminaries

This section introduced some basic concepts of surface curves and definitions of some special types of such curves in Euclidean 3-space $E^3$.

Let $M$ be an oriented surface in 3-dimensional Euclidean space $E^3$, and let $\alpha(s): I \subset IR \to M$ be a unit speed surface curve on $M$ and $s$ be the arc-length parameter of $\alpha$. If we denote the Frenet frame of $\alpha$ by $\{T, N, B\}$, then the Frenet equations of $\alpha$ are given by

$$T' = \kappa N, \ N' = -\kappa T + \tau B, \ B' = -\tau N,$$

where $X' = dX/ds$, $X \in \{T, N, B\}$; $\kappa(s)$ is curvature(or first curvature function), $\tau(s)$ is torsion(or second curvature function) and $T$, $N$ and $B$ are unit tangent vector, principal normal vector and binormal vector of $\alpha$, respectively.

On the other hand, since the curve $\alpha$ is a surface curve, it has another orthonormal frame called Darboux frame and this frame is denoted by $\{T, V, U\}$, where $T$ is unit tangent of the curve, $U$ is unit surface normal of $M$ along curve $\alpha$ and $V$ is a unit vector defined by $V = U \times T$. Using the fact that the unit tangent $T$ is common in both Frenet frame and Darboux frame, the relation between these frames can be given as follows

$$\begin{bmatrix} T \\ V \\ U \end{bmatrix} = \begin{bmatrix} 1 & 0 & 0 \\ 0 & \cos\varphi & \sin\varphi \\ 0 & -\sin\varphi & \cos\varphi \end{bmatrix} \begin{bmatrix} T \\ N \\ B \end{bmatrix},$$

where $\varphi$ is the angle between the vectors $V$ and $N$. The Darboux equations of $\alpha$ are given by

$$\begin{bmatrix} T' \\ V' \\ U' \end{bmatrix} = \begin{bmatrix} 0 & k_g & k_n \\ -k_g & 0 & \tau_g \\ -k_n & -\tau_g & 0 \end{bmatrix} \begin{bmatrix} T \\ V \\ U \end{bmatrix}, \tag{1}$$

where $k_n$, $k_g$ and $\tau_g$ are called normal curvature, geodesic curvature and geodesic torsion of $\alpha$, respectively[17]. The relations between these curvatures and $\kappa$, $\tau$ are given as follows

$$k_g = \kappa \cos\varphi, \ k_n = \kappa \sin\varphi, \ \tau_g = \tau + \frac{d\varphi}{ds}.$$

In the differential geometry of surfaces, for a curve $\alpha$ lying on a surface $M$, the followings are well-known

    **i)** $\alpha$ is a geodesic curve $\Leftrightarrow k_g = 0$,

    **ii)** $\alpha$ is an asymptotic line $\Leftrightarrow k_n = 0$,

    **iii)** $\alpha$ is a principal line $\Leftrightarrow \tau_g = 0$.

**Definition 2.1([16,19]).** Let $\gamma$ be a unit speed curve on an oriented surface $M$ and $\{T,V,U\}$ be the Darboux frame along $\gamma$. The curve $\gamma$ is called a relatively normal-slant helix(or $V$-strip slant helix) if the vector field $V$ of $\gamma$ makes a constant angle with a fixed direction, i.e. there exists a fixed unit vector $d$ and a constant angle $\theta$ such that $\langle V,d \rangle = \cos\theta$.

**Theorem 2.2([16]).** *A unit-speed curve $\gamma$ on a surface with $(\tau_g, k_g) \neq (0,0)$ is a relatively normal-slant helix if and only if*

$$f(s) = \left( \frac{1}{(k_g^2 + \tau_g^2)^{3/2}} \left( \tau_g' k_g - k_g' \tau_g - k_n(k_g^2 + \tau_g^2) \right) \right)(s)$$

*is a constant function.*

**Definition 2.3([10,19]).** Let $\beta$ be a unit speed curve on an oriented surface $M$ and $\{T,V,U\}$ be the Darboux frame along $\beta$. The curve $\beta$ is called an isophote curve(or $U$-strip slant helix) if the vector field $U$ of $\beta$ makes a constant angle with a fixed direction, i.e. there exists a fixed unit vector $l$ and a constant angle $\sigma$ such that $\langle U,l \rangle = \cos\sigma$.

**Theorem 2.4([10]).** *A unit-speed curve $\beta$ on a surface is an isophote curve if and only if*

$$\cot\sigma = \pm \left( \frac{k_n^2}{(k_n^2 + \tau_g^2)^{3/2}} \left( \frac{\tau_g}{k_n} \right)' + \frac{k_g}{(k_n^2 + \tau_g^2)^{1/2}} \right)(s)$$

*is a constant function, where $k_n \neq 0$.*

**Definition 2.5([12]).** Let $\beta$ be a unit speed curve on an oriented surface $M$ and $\{T,V,U\}$ be the Darboux frame along $\beta$. The vector fields $D_o(s)$, $D_n(s)$, $D_r(s)$ along $\beta$ defined by

$$D_o(s) = \tau_g(s)T(s) - k_n(s)V(s),$$
$$D_n(s) = -k_n(s)V(s) + k_g(s)U(s),$$
$$D_r(s) = \tau_g(s)T(s) + k_g(s)U(s),$$

are called osculating Darboux vector field, normal Darboux vector field and rectifying Darboux vector field along $\beta$, respectively.

Now, considering last definition, we can define a new type of surface curves as follows:

**Definition 2.6( $D_i$-Darboux slant helix).** Let $\beta$ be a unit speed curve on an oriented surface $M$ with Darboux vector fields $D_o(s)$, $D_n(s)$, $D_r(s)$. Then, $\beta$ is called $D_i$-Darboux slant helix if there exists a constant angle $\rho$ between the Darboux vector filed $D_i$ (or equivalently, unit Darboux vector field $\bar{D}_i(s) = D_i(s)/\|D_i(s)\|$) and a fixed(constant) unit direction $\zeta$, i.e., $\langle \bar{D}_i(s), \zeta \rangle = \cos\rho$ is constant, where $i \in \{o, n, r\}$.

## 3. Helices associated to surface curves

In this section, we introduce helices which are associated to special surface curves such as helical curves, relatively normal-slant helix and isophote curve. We give parametrizations of some special associated helices by means of Darboux frame of surface curve.

Let consider oriented surface $M$ and unit speed surface curve $\alpha : I \to M$ with arc-length parameter $s$ and let $\{T, V, U\}$ denotes the Darboux frame of $\alpha$. Consider the curve $\gamma : J \to E^3$ given by the parametrization

$$\gamma(s) = \alpha(s) + y_1(s)T(s) + y_2(s)V(s) + y_3(s)U(s), \tag{2}$$

where $y_i = y_i(s)$; $(1 \le i \le 3)$ are smooth functions of arc-length parameter $s$. The curve $\gamma$ is called an associated curve of surface curve $\alpha$ or surface curve connected(SCC) associated curve. From this definition, it is clear that it is not necessary for the curve $\gamma$ to lie on surface $M$, it can be lying on it, or not. Furthermore, if the associated curve $\gamma$ is a general helix, then it is called associated helix connected to surface curve $\alpha$ or SCC-associated helix.

Differentiating (2) with respect to $s$ and considering Darboux formulas (1), we get

$$\gamma'(s) = R_1(s)T(s) + R_2(s)V(s) + R_3(s)U(s), \tag{3}$$

where $R_i = R_i(s)$; $(1 \le i \le 3)$ are smooth functions of $s$ defined by

$$R_1 = y_1' - k_g y_2 - k_n y_3 + 1, \quad R_2 = y_2' + k_g y_1 - \tau_g y_3, \quad R_3 = y_3' + k_n y_1 + \tau_g y_2. \tag{4}$$

Considering (4), we can introduce special helices associated to special surface curves as follows:

### 3.1. Helices associated to helical curves on a surface

Let assume that $\alpha$ is any arbitrary curve on an oriented surface $M$ and $\gamma$ is an associated curve of $\alpha$. Let investigate the special case that tangent vector $\gamma'$ of $\gamma$ is linearly dependent with unit tangent vector $T$ of $\alpha$. For this, we have $R_1 \ne 0$, $R_2 = R_3 = 0$ and from (3), we can write $\gamma'(s) = R_1(s)T(s)$ and for the arc-length parameter $s_\gamma$ of curve $\gamma$, we have $ds_\gamma = \pm R_1(s)ds$. Then, we compute the Frenet vectors of associated curve $\gamma$ as follows:

$$\begin{cases} T_\gamma(s) = \pm T(s), \quad N_\gamma(s) = \dfrac{1}{\sqrt{k_g^2(s) + k_n^2(s)}}\left(k_g(s)V(s) + k_n(s)U(s)\right), \\ B_\gamma(s) = \dfrac{\pm 1}{\sqrt{k_g^2(s) + k_n^2(s)}}\left(-k_n(s)V(s) + k_g(s)U(s)\right) = \pm \dfrac{D_n(s)}{\|D_n(s)\|}, \end{cases} \tag{5}$$

It is well-known that if $\gamma$ is a general helix, then the unit tangent vector $T_\gamma$ makes a constant angle $\theta$ with a fixed direction $\xi$ in the form $\xi = (\cos\theta)T_\gamma(s) + (\sin\theta)B_\gamma(s)$. Then, from Definition 2.6 and equation (5), the following theorem can be given:

**Theorem 3.1.** *Let $\gamma$ be an associated curve such that tangent vector $\gamma'$ of $\gamma$ is linearly dependent with unit tangent vector $T$ of surface curve $\alpha$ with $(k_g(s), k_n(s)) \ne (0,0)$. The followings are equivalent,*
  *(i) $\gamma$ is a helix.*

*(ii)* $\alpha$ *is a helical curve on* $M$.
*(iii)* $\alpha$ *is a* $D_n$*-Darboux slant helix on* $M$.

**Corollary 3.2.** *The helix* $\gamma$ *associated to helical curve* $\alpha$ *can be named by: Helical curve-connected associated helix or HCC-associated helix.*

From (4), the following system is obtained,
$$y_1' - k_g y_2 - k_n y_3 + 1 \neq 0, \ y_2' + k_g y_1 - \tau_g y_3 = 0, \ y_3' + k_n y_1 + \tau_g y_2 = 0. \tag{6}$$
Now, to obtain the parametric form of associated curve $\gamma$, we investigate following special cases:

**Case (1):** $y_1 = 0$. Then, system (6) becomes
$$k_g y_2 - k_n y_3 + 1 \neq 0, \ y_2' - \tau_g y_3 = 0, \ y_3' + \tau_g y_2 = 0,$$
and by solving this system, we get
$$y_1 = 0, \ y_2 = \sin \int \tau_g(s)ds, \ y_3 = \cos \int \tau_g(s)ds.$$
Writing these results in (2) gives that the parametric representation of SCC-associated curve $\gamma$ is obtained as
$$\gamma(s) = \alpha(s) + \left[\sin\left(\int \tau_g(s)ds\right)\right]V(s) + \left[\cos\left(\int \tau_g(s)ds\right)\right]U(s), \tag{7}$$
and we can give the followings:

**Theorem 3.3.** *The associated curve* $\gamma$ *given in (7) is a helix curve(or general helix) if and only if* $\alpha$ *is a helical curve on* $M$.
**Corollary 3.4.** *The associated curve (7) can be named by: Helical curve-connected associated helix of type 1 or HCC-associated helix of type 1.*
**Corollary 3.5.** *The helical curve* $\alpha$ *is a principal line if and only if HCC-associated helix (7) has the parametrization* $\gamma(s) = \alpha(s) + U(s)$.

**Case (2):** $y_2 = 0$. For this case, from (6) we have
$$y_1' - k_n y_3 + 1 \neq 0, \ k_g y_1 - \tau_g y_3 = 0, \ y_3' + k_n y_1 = 0. \tag{8}$$
In this system, if $k_g = 0$, we have $\tau_g y_3 = 0$, $y_3' + k_n y_1 = 0$. For this new case, if $\tau_g = 0$, we cannot determinate the functions $y_1(s), y_3(s)$. Taking $\tau_g \neq 0$, we get $y_3 = 0$, $k_n y_1 = 0$. Now, if $k_n = 0$, we cannot determine the function $y_1(s)$ and if $k_n \neq 0$, we get $y_1 = 0$ and $\gamma(s) = \alpha(s)$. So, we investigate the solution for $k_g \neq 0$. Then, From the second equation of (8), we get $y_1 = (\tau_g / k_g)y_3$. By writing last equality in the third equation of (8), it follows $y_3' + \dfrac{k_n \tau_g}{k_g} y_3 = 0$. The solution of this differential equation is $y_3 = c_1 \exp\left(-\int \dfrac{k_n \tau_g}{k_g} ds\right)$, where $c_1$ is a non-zero real constant. So, the solution of system (6) is
$$y_1 = c_1 \frac{\tau_g}{k_g} \exp\left(-\int \frac{k_n \tau_g}{k_g} ds\right), \ y_2 = 0, \ y_3 = c_1 \exp\left(-\int \frac{k_n \tau_g}{k_g} ds\right).$$
Then, from (2), the parametric representation of associated curve $\gamma$ is obtained as

$$\gamma(s) = \alpha(s) + c_1 \exp\left(-\int \frac{k_n \tau_g}{k_g} ds\right)\left(\frac{\tau_g}{k_g} T(s) + U(s)\right), \qquad (9)$$

and we can write the followings:

**Theorem 3.6.** *The associated curve $\gamma$ given in (9) is a helix curve(or general helix) if and only if $\alpha$ is a helical curve with $k_g \neq 0$ on $M$.*

**Corollary 3.7.** *The associated curve (9) can be named by: Helical curve-connected associated helix of type 2 or HCC-associated helix of type 2.*

**Corollary 3.8.** *(i) The helical curve $\alpha$ with $k_g \neq 0$ is an asymptotic curve if and only if HCC-associated helix (9) has the parametrization $\gamma(s) = \alpha(s) + c_1 \left(\frac{\tau_g}{k_g} T(s) + U(s)\right)$.*

*(ii) The helical curve $\alpha$ with $k_g \neq 0$ is a principal line if and only if HCC-associated helix (9) has the parametrization $\gamma(s) = \alpha(s) + c_1 U(s)$.*

**Case (3):** $y_3 = 0$. In this case, system (6) reduces to
$$y_1' - k_g y_2 + 1 \neq 0, \quad y_2' + k_g y_1 = 0, \quad k_n y_1 + \tau_g y_2 = 0. \qquad (10)$$
In this system, if $k_n = 0$, we have $y_2' + k_g y_1 = 0$, $\tau_g y_2 = 0$. For this new case, if $\tau_g = 0$, we cannot determinate the functions $y_1(s)$, $y_2(s)$. Taking $\tau_g \neq 0$, we get $k_g y_1 = 0$, $y_2 = 0$. Now, if $k_g = 0$, we cannot determinate the function $y_1(s)$ and if $k_g \neq 0$, we get $y_1 = 0$ and $\gamma(s) = \alpha(s)$. So, we investigate the solution for $k_n \neq 0$. From the third equation of (10), we get $y_1 = -(\tau_g / k_n) y_2$. By writing last equality into the second equation of (10), it follows $y_2' - \frac{k_g \tau_g}{k_n} y_2 = 0$. The solution of this differential equation is $y_2 = c_2 \exp\left(\int \frac{k_g \tau_g}{k_n} ds\right)$, where $c_2$ is a non-zero real constant. So, the solution of system (6) is
$$y_1 = -c_2 \frac{\tau_g}{k_n} \exp\left(\int \frac{k_g \tau_g}{k_n} ds\right), \quad y_2 = c_2 \exp\left(\int \frac{k_g \tau_g}{k_n} ds\right), \quad y_3 = 0.$$
Then, from (2), the parametric representation of associated curve $\gamma$ is obtained as
$$\gamma(s) = \alpha(s) - c_2 \exp\left(\int \frac{k_g \tau_g}{k_n} ds\right)\left(\frac{\tau_g}{k_n} T(s) - V(s)\right), \qquad (11)$$
and we can write the followings:

**Theorem 3.9.** *The associated curve $\gamma$ given in (11) is a helix curve(or general helix) if and only if $\alpha$ is a helical curve with $k_n \neq 0$ on $M$.*

**Corollary 3.10.** *The associated curve (11) can be named by: Helical curve-connected associated helix of type 3 or HCC-associated helix of type 3.*

**Corollary 3.11.** *(i) The helical curve $\alpha$ with $k_n \neq 0$ is a geodesic curve if and only if HCC-associated helix (11) has the parametrization $\gamma(s) = \alpha(s) - c_2 \left(\frac{\tau_g}{k_n} T(s) - V(s)\right)$.*

*(ii) The helical curve $\alpha$ with $k_n \neq 0$ is a principal line if and only if HCC-associated helix (11) has the parametrization $\gamma(s) = \alpha(s) + c_2 V(s)$.*

## 3.2. Helices associated to relatively normal-slant helices

In this section, we will consider special helices associated to relatively normal-slant helices. For this, first assume that the tangent vector $\gamma'$ of $\gamma$ is linearly dependent with unit vector field $V$ of surface curve $\alpha$. From (3), we can write $\gamma'(s) = R_2(s)V(s)$ and for the arc-length parameter $s_\gamma$ of curve $\gamma$, we have $ds_\gamma = \pm R_2(s)ds$. Then, after some computations, Frenet vectors $T_\gamma(s)$, $N_\gamma(s)$, $B_\gamma(s)$ of associated curve $\gamma$ are calculated as follows:

$$\begin{cases} T_\gamma(s) = \pm V(s), \ N_\gamma(s) = \dfrac{1}{\sqrt{k_g^2(s) + \tau_g^2(s)}} \left( -k_g(s)T(s) + \tau_g(s)U(s) \right), \\ B_\gamma(s) = \dfrac{\pm 1}{\sqrt{k_g^2(s) + \tau_g^2(s)}} \left( \tau_g(s)T(s) + k_g(s)U(s) \right) = \pm \dfrac{D_r(s)}{\|D_r(s)\|}. \end{cases} \quad (12)$$

If $\gamma$ is a general helix, it is well-known that both unit tangent vector $T_\gamma$ and binormal vector $B_\gamma$ make a constant angle $\theta$ with the same fixed direction. Then, from Definition 2.6 and equation (12), the following theorem can be given:

**Theorem 3.12.** *Let $\gamma$ be an associated curve such that tangent vector $\gamma'$ of $\gamma$ is linearly dependent with unit vector field $V$ of surface curve $\alpha$ with $(k_g(s), \tau_g(s)) \neq (0,0)$. The followings are equivalent,*
  *(i) $\gamma$ is a helix.*
  *(ii) $\alpha$ is a relatively normal-slant helix.*
  *(iii) $\alpha$ is a $D_r$-Darboux slant helix.*

**Corollary 3.13.** *The helix $\gamma$ associated to relatively normal-slant helix $\alpha$ can be named by: Relatively normal-slant helix-connected associated helix or RNS-HC-associated helix.*

From (4), the following system is obtained
$$y_2' + k_g y_1 - \tau_g y_3 \neq 0, \ y_1' - k_g y_2 - k_n y_3 + 1 = 0, \ y_3' + k_n y_1 + \tau_g y_2 = 0. \quad (13)$$

Now, to obtain the parametric form of associated curve $\gamma$, we investigate following special cases:

**Case (i):** $y_1 = 0$. Then, system (13) becomes
$$y_2' - \tau_g y_3 \neq 0, \ -k_g y_2 - k_n y_3 + 1 = 0, \ y_3' + \tau_g y_2 = 0. \quad (14)$$

If $k_g = 0, k_n \neq 0, \tau_g \neq 0$, the solution of system (14) is $y_3 = 1/k_n$, $y_2 = k_n'/\tau_g k_n^2$. Then, the parametric form of $\gamma$ is

$$\gamma(s) = \alpha(s) + \left( \dfrac{k_n'}{k_n^2 \tau_g} \right)(s)V(s) + \dfrac{1}{k_n(s)}U(s) \quad (15)$$

If $k_g \neq 0$, from the second equation of (14), we have $y_2 = -(k_n/k_g)y_3 + (1/k_g)$. Writing last equality in the third equation of (14) gives $y_3' - (k_n \tau_g / k_g)y_3 + (\tau_g / k_g) = 0$ and by solving this differential equation, the solution of system (14) is

$$y_2 = \frac{k_n}{k_g}\exp\left(\int\frac{k_n\tau_g}{k_g}ds\right)\left[\int\exp\left(-\int\frac{k_n\tau_g}{k_g}ds\right)\frac{\tau_g}{k_g}ds - c_3\right] + \frac{1}{k_g},$$

$$y_3 = \exp\left(\int\frac{k_n\tau_g}{k_g}ds\right)\left[c_3 - \int\exp\left(-\int\frac{k_n\tau_g}{k_g}ds\right)\frac{\tau_g}{k_g}ds\right],$$

where $c_3$ is a real constant. By writing these results in (2), the parametric representation of SCC-associated curve $\gamma$ is obtained as

$$\gamma(s) = \alpha(s) + \left(\frac{k_n}{k_g}\exp\left(\int\frac{k_n\tau_g}{k_g}ds\right)\left[\int\exp\left(-\int\frac{k_n\tau_g}{k_g}ds\right)\frac{\tau_g}{k_g}ds - c_3\right] + \frac{1}{k_g}\right)V(s)$$

$$+ \exp\left(\int\frac{k_n\tau_g}{k_g}ds\right)\left[c_3 - \int\exp\left(-\int\frac{k_n\tau_g}{k_g}ds\right)\frac{\tau_g}{k_g}ds\right]U(s)$$

(16)

and we can give the followings:

**Theorem 3.14.** *The associated curve $\gamma$ given in (15) and (16) is a general helix if and only if $\alpha$ is a relatively normal-slant helix on $M$.*

**Corollary 3.15.** *The associated curves (15) and (16) can be named by: Relatively normal-slant helix-connected associated helix of type 1 or RNS-HC-associated helix of type 1.*

**Corollary 3.16.** *(i) Relatively normal-slant helix $\alpha$ with $k_g \neq 0$ is an asymptotic curve if and only if RNS-HC-associated helix (16) has the parametrization*

$$\gamma(s) = \alpha(s) + \frac{1}{k_g}V(s) + \left[c_3 - \int\frac{\tau_g}{k_g}ds\right]U(s).$$

*(ii) Relatively normal-slant helix $\alpha$ with $k_g \neq 0$ is a principal line if and only if RNS-HC-associated helix (16) has the parametrization* $\gamma(s) = \alpha(s) + \left(\frac{1-c_3 k_n}{k_g}\right)V(s) + c_3 U(s).$

**Case (ii):** $y_2 = 0$. For this case, from (13) we have

$$k_g y_1 - \tau_g y_3 \neq 0,\quad y_1' - k_n y_3 + 1 = 0,\quad y_3' + k_n y_1 = 0. \tag{17}$$

In this system, if $k_n = 0$, we have $y_1 = c_4 - s$, $y_3 = c_5$, where $c_4, c_5$ are integration constants. Then, the associated curve $\gamma$ has the form

$$\gamma(s) = \alpha(s) + (c_4 - s)T(s) + c_5 U(s). \tag{18}$$

If $k_n \neq 0$, differentiating third equation in (17) gives

$$y_1' = \frac{k_n'}{k_n^2}y_3' - \frac{1}{k_n}y_3''. \tag{19}$$

By writing (19) into the second equation of (17), we have

$$y_3'' - \frac{k_n'}{k_n}y_3' + k_n^2 y_3 - k_n = 0. \tag{20}$$

By changing the variable as $\theta = \int k_n(s)ds$, the homogeneous part of (20) becomes

$$\frac{d^2 y_3}{d\theta^2} + y_3 = 0,$$

and it has the solution

$$y_{3h}(s) = c_6 \cos\left(\int k_n(s)ds\right) + c_7 \sin\left(\int k_n(s)ds\right), \qquad (21)$$

where $c_6, c_7$ are real constants. Now, considering (21) and by using the variation of parameters method, the solution of (20) is

$$y_3(s) = c_6 \cos\left(\int k_n(s)ds\right) + c_7 \sin\left(\int k_n(s)ds\right) - \cos\left(\int k_n(s)ds\right)\int\left[\sin\left(\int k_n(s)ds\right)\right]ds \\ + \sin\left(\int k_n(s)ds\right)\int\left[\cos\left(\int k_n(s)ds\right)\right]ds \qquad (22)$$

Differentiating (22) and writing the result in the third equation in (17), we have

$$y_1(s) = -\sin\left(\int k_n(s)ds\right)\left[\int \sin\left(\int k_n(s)ds\right)ds - c_6\right] \\ - \cos\left(\int k_n(s)ds\right)\left[\int \cos\left(\int k_n(s)ds\right)ds + c_7\right] \qquad (23)$$

Then, the parametric form of associated curve $\gamma$ is given by

$$\gamma(s) = \alpha(s) + y_1(s)T(s) + y_3(s)U(s), \qquad (24)$$

where $y_1(s)$ and $y_3(s)$ are given in (23) and (22), respectively. Now, we can write the followings:

**Theorem 3.17.** *The associated curve $\gamma$ given in (18) and (24) is a general helix if and only if $\alpha$ is a relatively normal-slant helix on $M$.*

**Corollary 3.18.** *The associated curves (18) and (24) can be named by: Relatively normal-slant helix-connected associated helix of type 2 or RNS-HC-associated helix of type 2.*

**Case (iii):** $y_3 = 0$. Then, system (13) reduces to

$$y_2' + k_g y_1 \neq 0, \ y_1' - k_g y_2 + 1 = 0, \ k_n y_1 + \tau_g y_2 = 0. \qquad (25)$$

In this system, if $\tau_g = 0$, $k_n \neq 0$, $k_g \neq 0$, we have, $y_1 = 0$, $y_2 = 1/k_g$ and the associated curve $\gamma$ has the parametric form

$$\gamma(s) = \alpha(s) + (1/k_g)(s)V(s). \qquad (26)$$

If $\tau_g \neq 0$, from the third equation of (25), we get $y_2 = -(k_n/\tau_g)y_1$. By writing last equality into the second equation of (25), it follows $y_1' + \dfrac{k_n k_g}{\tau_g}y_1 + 1 = 0$. The solution of this differential equation is

$$y_1 = \left(\exp\left(-\int \frac{k_n k_g}{\tau_g}ds\right)\right)\left(c_8 - \int \exp\left(\int \frac{k_n k_g}{\tau_g}ds\right)ds\right), \qquad (27)$$

where $c_8$ is a real constant. Then, from the third equation of (25), we have

$$y_2 = -\frac{k_n}{\tau_g}\left(\exp\left(-\int \frac{k_n k_g}{\tau_g}ds\right)\right)\left(c_8 - \int \exp\left(\int \frac{k_n k_g}{\tau_g}ds\right)ds\right). \qquad (28)$$

Then, from (2), the parametric representation of associated curve $\gamma$ is obtained as

$$\gamma(s) = \alpha(s) + y_1(s)T(s) + y_2(s)V(s), \qquad (29)$$

where $y_1(s)$ and $y_2(s)$ are defined in (27) and (28), respectively. Now, we can write the followings:

**Theorem 3.19.** *The associated curve $\gamma$ given in (26) and (29) is a helix curve(or general helix) if and only if $\alpha$ is a relatively normal-slant helix on $M$.*

**Corollary 3.20.** *The associated curves (26) and (29) can be named by: Relatively normal-slant helix-connected associated helix of type 3 or RNS-HC-associated helix of type 3.*

**Corollary 3.21.** *(i) Relatively normal-slant helix $\alpha$ with $\tau_g \neq 0$ is a geodesic curve if and only if RNS-HC-associated helix (29) has the parametrization*

$$\gamma(s) = \alpha(s) + (c_8 - s)\big(T(s) - (k_n/\tau_g)(s)V(s)\big).$$

*(ii) Relatively normal-slant helix $\alpha$ with $\tau_g \neq 0$ is an asymptotic curve if and only if RNS-HC-associated helix (29) has the parametrization $\gamma(s) = \alpha(s) + (c_8 - s)T(s)$.*

### 3.3. Helices associated to isophote curves on a surface

Let assume that $\alpha$ is any arbitrary curve on an oriented surface $M$ and $\gamma$ is an associated curve of $\alpha$. Let investigate the special case that tangent vector $\gamma'$ of $\gamma$ is linearly dependent with unit vector field $U$ of $\alpha$. For this, we have $R_1 = R_2 = 0$, $R_3 \neq 0$ and from (3), we can write $\gamma'(s) = R_3(s)U(s)$. Then, the Frenet vectors of associated curve $\gamma$ are calculated as follows:

$$\begin{cases} T_\gamma(s) = \pm U(s), \quad N_\gamma(s) = \dfrac{1}{\sqrt{k_n^2(s) + \tau_g^2(s)}}\big(-k_n(s)T(s) - \tau_g(s)V(s)\big), \\ B_\gamma(s) = \dfrac{\pm 1}{\sqrt{k_n^2(s) + \tau_g^2(s)}}\big(\tau_g(s)T(s) - k_n(s)V(s)\big) = \pm \dfrac{D_o(s)}{\|D_o(s)\|}, \end{cases} \quad (30)$$

and from Definition 2.6 and equation (30), the followings can be given:

**Theorem 3.22.** *Let $\gamma$ be an associated curve such that tangent vector $\gamma'$ of $\gamma$ is linearly dependent with unit vector field $U$ of surface curve $\alpha$ with $(k_n(s), \tau_g(s)) \neq (0,0)$. The followings are equivalent,*
   *(i) $\gamma$ is a helix.*
   *(ii) $\alpha$ is an isophote curve.*
   *(iii) $\alpha$ is a $D_o$-Darboux slant helix.*

**Corollary 3.23.** *The helix $\gamma$ associated to isophote curve $\alpha$ can be named by: Isophote curve connected associated helix or ICC-associated helix.*

From (4), the following system is obtained
$$y_3' + k_n y_1 + \tau_g y_2 \neq 0, \quad y_1' - k_g y_2 - k_n y_3 + 1 = 0, \quad y_2' + k_g y_1 - \tau_g y_3 = 0, \quad (31)$$
and to determinate the functions $y_i(s)$, $(1 \leq i \leq 3)$, we investigate following special cases:

**Case (A):** $y_1 = 0$. Then, system (31) becomes
$$y_3' + \tau_g y_2 \neq 0, \quad -k_g y_2 - k_n y_3 + 1 = 0, \quad y_2' - \tau_g y_3 = 0. \quad (32)$$
In this system, if $k_n = 0$ and $k_g \neq 0$, $\tau_g \neq 0$, we have $y_2 = -1/k_g$, $y_3 = k_g'/k_g^2 \tau_g$. Then, the parametric form of associated curve $\gamma$ is obtained as

$$\gamma(s) = \alpha(s) - \dfrac{1}{k_g(s)}V(s) + \left(\dfrac{k_g'}{k_g^2 \tau_g}\right)(s)U(s). \quad (33)$$

If $k_n \neq 0$, from the second equation in (32), we get $y_3 = -(k_g/k_n)y_2 + (1/k_n)$. By writing last equality in the third equation in (32), it follows $y_2' + \frac{k_g \tau_g}{k_n} y_2 - \frac{\tau_g}{k_n} = 0$. The solution of this differential equation is

$$y_2 = \exp\left(-\int \frac{k_g \tau_g}{k_n} ds\right)\left[\int \exp\left(\int \frac{k_g \tau_g}{k_n} ds\right) \frac{\tau_g}{k_n} ds + c_8\right], \tag{34}$$

where $c_8$ is a real constant. By writing (34) into the second equality of (32), we have

$$y_3 = -\frac{k_g}{k_n} \exp\left(-\int \frac{k_g \tau_g}{k_n} ds\right)\left[\int \exp\left(\int \frac{k_g \tau_g}{k_n} ds\right) \frac{\tau_g}{k_n} ds + c_8\right] + \frac{1}{k_n}. \tag{35}$$

Then, from (2), the parametric representation of associated curve $\gamma$ is obtained as

$$\gamma(s) = \alpha(s) + y_2(s)V(s) + y_3(s)U(s), \tag{36}$$

where $y_2(s)$ and $y_3(s)$ are given in (34) and (35), respectively, and we can write the followings:

**Theorem 3.24.** *The associated curve $\gamma$ given in (33) and (36) is a helix curve(or general helix) if and only if $\alpha$ is an isophote curve on $M$.*

**Corollary 3.25.** *The associated curves (33) and (36) can be named by: Isophote curve-connected associated helix of type 1 or ICC-associated helix of type 1.*

**Corollary 3.26.** *(i) The isophote curve $\alpha$ with $k_n \neq 0$ is a geodesic curve if and only if ICC-associated helix (36) has the parametrization $\gamma(s) = \alpha(s) + \left(\int \frac{\tau_g}{k_n} ds + c_8\right)V(s) + \frac{1}{k_n(s)} U(s)$.*

*(ii) The isophote curve $\alpha$ with $k_n \neq 0$ is a principal line if and only if ICC-associated helix (36) has the parametrization $\gamma(s) = \alpha(s) + c_8 V(s) + \frac{1 - c_8 k_g(s)}{k_n(s)} U(s)$.*

**Case (B):** $y_2 = 0$. In this case, system (31) reduces to

$$y_3' + k_n y_1 \neq 0, \quad y_1' - k_n y_3 + 1 = 0, \quad k_g y_1 - \tau_g y_3 = 0. \tag{37}$$

In this system, if $\tau_g = 0$, we have $y_1' - k_n y_3 + 1 = 0$, $k_g y_1 = 0$. For this new case, if $k_g = 0$, we cannot determinate the functions $y_1(s), y_3(s)$. Taking $k_g \neq 0$, $k_n \neq 0$, we get $y_3 = -1/k_n$, $y_1 = 0$. Then, the parametric form of associated curve $\gamma$ is

$$\gamma(s) = \alpha(s) - \frac{1}{k_n(s)} U(s). \tag{38}$$

If $\tau_g \neq 0$, from the third equation of (37), we get $y_3 = (k_g/\tau_g)y_1$. By writing last equality into the second equation of (37), it follows $y_1' - \frac{k_g k_n}{\tau_g} y_1 + 1 = 0$. The solution of this differential equation is

$$y_1 = \exp\left(\int \frac{k_g k_n}{\tau_g} ds\right)\left[c_9 - \int \exp\left(-\int \frac{k_g k_n}{\tau_g} ds\right) ds\right], \tag{39}$$

where $c_9$ is a real constant. Writing (39) into the third equation of (37) gives

$$y_3 = \frac{k_g}{\tau_g} \exp\left(\int \frac{k_g k_n}{\tau_g} ds\right) \left[c_9 - \int \exp\left(-\int \frac{k_g k_n}{\tau_g} ds\right) ds\right]. \tag{40}$$

From (2), the parametric representation of associated curve $\gamma$ is obtained as

$$\gamma(s) = \alpha(s) + y_1(s)T(s) + y_3(s)V(s), \tag{41}$$

where $y_1(s)$ and $y_3(s)$ are given in (39) and (40), respectively. So, we can write the followings:

**Theorem 3.27.** *The associated curve $\gamma$ given in (38) and (41) is a helix curve(or general helix) if and only if $\alpha$ is an isophote curve on $M$.*

**Corollary 3.28.** *The associated curves (38) and (41) can be named by: Isophote curve-connected associated helix of type 2 or ICC-associated helix of type 2.*

**Corollary 3.29.** *(i) The isophote curve $\alpha$ with $\tau_g \neq 0$ is a geodesic curve if and only if ICC-associated helix (41) has the parametrization $\gamma(s) = \alpha(s) + (c_9 - s)T(s)$.*

*(ii) The isophote curve $\alpha$ with $\tau_g \neq 0$ is an asymptotic curve if and only if ICC-associated helix (41) has the parametrization $\gamma(s) = \alpha(s) + (c_9 - s)\left(T(s) + \frac{k_g(s)}{\tau_g(s)}U(s)\right)$.*

**Case (C):** $y_3 = 0$. For this case, from (31) we have

$$k_n y_1 + \tau_g y_2 \neq 0, \quad y_1' - k_g y_2 + 1 = 0, \quad y_2' + k_g y_1 = 0. \tag{42}$$

In this system, if $k_g = 0$, we have $y_1 = c_{10} - s$, $y_2 = c_{11}$, where $c_{10}$ and $c_{11}$ are integration constants. Then, the associated curve $\gamma$ has the form

$$\gamma(s) = \alpha(s) + (c_{10} - s)T(s) + c_{11}V(s). \tag{43}$$

If $k_g \neq 0$, differentiating third equation in (42) gives

$$y_1' = \frac{k_g'}{k_g^2} y_2' - \frac{1}{k_g} y_2''. \tag{44}$$

By writing (44) into the second equation of (43), we have

$$y_2'' - \frac{k_g'}{k_g} y_2' + k_g^2 y_2 - k_g = 0. \tag{45}$$

By changing the variable as $\varsigma = \int k_g(s) ds$, the homogeneous part of (45) becomes

$$\frac{d^2 y_2}{d\varsigma^2} + y_2 = 0,$$

and it has the solution

$$y_{2h}(s) = c_{12} \cos\left(\int k_g(s) ds\right) + c_{13} \sin\left(\int k_g(s) ds\right), \tag{46}$$

where $c_{12}, c_{13}$ are real constants. Now, considering (46) and by using the variation of parameters method, the solution of (45) is

$$y_2(s) = c_{12} \cos\left(\int k_g(s) ds\right) + c_{13} \sin\left(\int k_g(s) ds\right) - \cos\left(\int k_g(s) ds\right) \int \left[\sin\left(\int k_g(s) ds\right)\right] ds \\ + \sin\left(\int k_g(s) ds\right) \int \left[\cos\left(\int k_g(s) ds\right)\right] ds \tag{47}$$

Differentiating (47) and writing the result in the third equation in (43), we have

$$y_1(s) = -\sin\left(\int k_g(s)ds\right)\left[\int \sin\left(\int k_g(s)ds\right)ds - c_{12}\right]$$
$$-\cos\left(\int k_g(s)ds\right)\left[\int \cos\left(\int k_g(s)ds\right)ds + c_{13}\right] \quad (48)$$

Then, the parametric form of associated curve $\gamma$ is given by

$$\gamma(s) = \alpha(s) + y_1(s)T(s) + y_2(s)V(s), \quad (49)$$

where $y_1(s)$ and $y_2(s)$ are given in (48) and (47), respectively. Now, we get the followings:

**Theorem 3.30.** *The associated curve $\gamma$ given in (43) and (49) is a general helix if and only if $\alpha$ is an isophote curve on $M$.*

**Corollary 3.31.** *The associated curves (43) and (49) can be named by: Isophote curve-connected associated helix of type 3 or ICC-associated helix of type 3.*

## 4. Examples

**Example 4.1.** Let consider the cylinder surface $M$ given by the parametrization $\varphi(u,v) = (\sin u, \cos u, v)$ and helical curve $\alpha : I \to M$ defined by the parametric form $\alpha(s) = \left(\sin\frac{s}{\sqrt{2}}, \cos\frac{s}{\sqrt{2}}, \frac{s}{\sqrt{2}}\right)$. The vectors of Darboux frame and curvatures of $\alpha$ are computed as follows,

$$T(s) = \frac{1}{\sqrt{2}}\left(\cos\frac{s}{\sqrt{2}}, -\sin\frac{s}{\sqrt{2}}, 1\right), \quad V(s) = \frac{1}{\sqrt{2}}\left(-\cos\frac{s}{\sqrt{2}}, \sin\frac{s}{\sqrt{2}}, 1\right),$$

$$U(s) = \left(-\sin\frac{s}{\sqrt{2}}, -\cos\frac{s}{\sqrt{2}}, 0\right); \quad k_g = 0, \; k_n = 1/2, \; \tau_g = -1/2,$$

i.e., $\alpha$ is a geodesic. Since $V(s)$ and $U(s)$ are circles, it is clear that this curve is also both a relatively normal-slant helix and isophote curve on $M$. Then, using the obtained results, some special helices associated to $\alpha$ are given in Fig. 1.

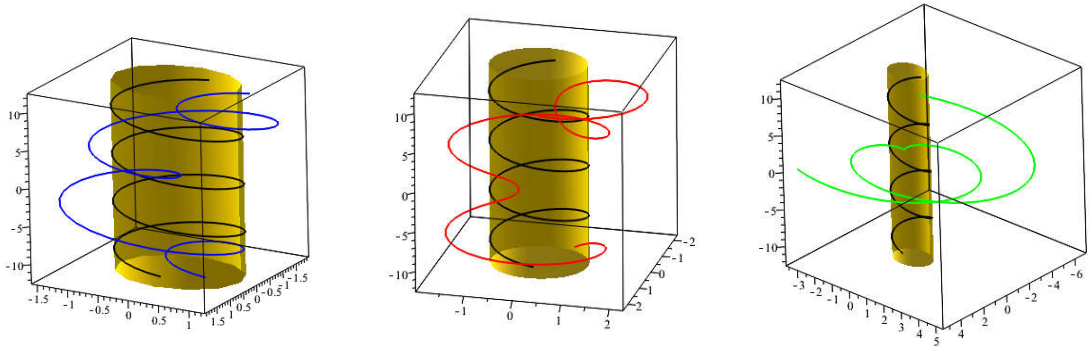

**Fig. 1.** Surface curve $\alpha$ (black), HCC-associated helix of type 1(blue), RNS-HC-associated helix of type 2(red) and ICC-associated helix of type 1(green), respectively.

**Example 4.2.** Let $M$ be a helicoid surface given by the parametrization $\varphi(u,v) = (v\cos u, v\sin u, u)$ and consider the helical curve $\alpha : I \to M$ defined by $\alpha(s) = \left(\cos\frac{s}{\sqrt{2}}, \sin\frac{s}{\sqrt{2}}, \frac{s}{\sqrt{2}}\right)$. For this curve, we have,

$$T(s) = \frac{1}{\sqrt{2}}\left(-\sin\frac{s}{\sqrt{2}}, \cos\frac{s}{\sqrt{2}}, 1\right), \quad V(s) = \left(\cos\frac{s}{\sqrt{2}}, \sin\frac{s}{\sqrt{2}}, 0\right),$$

$$U(s) = \frac{1}{\sqrt{2}}\left(-\sin\frac{s}{\sqrt{2}}, \cos\frac{s}{\sqrt{2}}, -1\right); \quad k_g = -1/2, \; k_n = 0, \; \tau_g = 1/2,$$

i.e., $\alpha$ is an asymptotic curve. Since $V(s)$ and $U(s)$ are circles, we have that this curve is also both a relatively normal-slant helix and isophote curve on $M$. Then, using the obtained results, some special helices associated to $\alpha$ are given in Fig. 2.

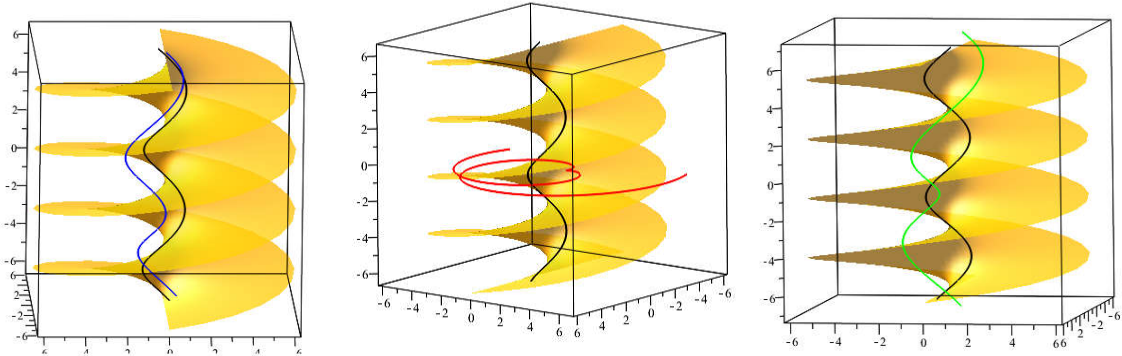

**Fig. 2.** Surface curve $\alpha$ (black), HCC-associated helix of type 1(blue), RNS-HC-associated helix of type 1(red) and ICC-associated helix of type 3(green), respectively.

**Conclusions**

By considering the special surface curves such as helical curves, relatively normal-slant helices and isophote curves, some special types of general helices are introduced. These new helices are called associated helices connected to surface curves and parametric representations of them are introduced by means of Darboux frame and surface curvatures of main surface curve.

**Author Declaration:** I wish to confirm that there are no known conflicts of interest associated with this publication and there has been no significant financial support for this work that could have influenced its outcome.